\documentclass[12pt]{article}

\usepackage[a4paper,margin=1in]{geometry}

\usepackage{amsmath,amsthm,amssymb,amsfonts}
\usepackage{mathrsfs}
\usepackage{extarrows}
\usepackage{mathtools} 

\usepackage{tikz}
\usetikzlibrary{positioning}
\usetikzlibrary{calc}
\usetikzlibrary{arrows.meta,decorations.pathreplacing,calligraphy}

\usepackage{xcolor}
\usepackage{enumerate}
\usepackage{float}
\usepackage{booktabs}

\usepackage[numbers,square,comma,sort&compress]{natbib}
\usepackage[hidelinks]{hyperref}

\newtheorem{thm}{Theorem}[section]
\newtheorem{lem}[thm]{Lemma}

\newtheorem{prop}[thm]{Proposition}

\theoremstyle{definition}
\newtheorem{defi}[thm]{Definition}

\newtheorem{rem}[thm]{Remark}

\newcommand{\defn}[1]{\emph{#1}}

\newcommand{\Del}{\Delta}

\newcommand{\C}{\mathbb{C}}
\newcommand{\R}{\mathbb{R}}

\newcommand{\bfx}{\mathbf{x}}

\newcommand{\cH}{\mathcal{H}}

\newcommand{\cP}{\mathcal{P}}
\newcommand{\cQ}{\mathcal{Q}}

\newcommand{\Ima}{{\rm Im}}

\newcommand{\wt}{{\rm wt}}

\newcommand{\what}{\widehat}
\newcommand{\wtil}{\widetilde}

\numberwithin{equation}{section}

\title{Real stable multivariate Delannoy polynomials and their Rayleigh differences}
\author{Xi Chen%
\thanks{
	Corresponding author.
	\newline\hspace*{3mm}
	{\it Email address:}\quad
	chenxi@dlut.edu.cn (X. Chen),
	ruiqimao@mail.dlut.edu.cn (R. Mao)
},
Ruiqi Mao
}
\date{\small
School of Mathematical Sciences, Dalian University of Technology, Dalian 116024, PR China}

\begin{document}
\maketitle

\begin{abstract}
The univariate Delannoy polynomials are known to have only real zeros. We construct a real stable multivariate refinement of these polynomials using weighted lattice paths and give both algebraic and combinatorial proofs of stability. More strongly, we show that every Rayleigh difference of the multivariate Delannoy polynomials is a perfect square and give an explicit factorization.
\end{abstract}

\section{Introduction}

A \defn{Delannoy path of size $n$} is a lattice path from $(0,0)$ to a point on $x+y=n$ using steps $H=(1,0)$, $V=(0,1)$, and $D=(1,1)$.
The Delannoy numbers $d(n,k)$ count the Delannoy paths from $(0,0)$ to $(n-k,k)$.
These numbers are named after Henri Delannoy (see \cite{BS05} for historical remarks).
The $n$-th \defn{Delannoy polynomial} is defined as
\[
d_n(x)=\sum_{k=0}^{n}d(n,k)x^k, \qquad n=0,1,2,\dots.
\]
Wang et al.~\cite{WZC19} proved that the zeros of the Delannoy polynomials are all real and distinct.
A natural multivariate extension of real-rootedness is real stability.
A polynomial $f(z_1,\dots,z_n)\in\C[z_1,\dots,z_n]$ is called \defn{stable} if $f(z_1,\dots,z_n)\neq 0$ whenever $\Ima (z_i)>0$ for $1\le i\le n$ or $f$ is identically zero.
If in addition $f$ has real coefficients, it is called \defn{real stable}.
It is known~\cite{Wag11} that a nonzero univariate polynomial is real stable if and only if it has only real zeros.

Stable multivariate refinements have been extensively studied for many classical combinatorial polynomials, 
including the Eulerian polynomials with various extensions~\cite{HV12,YL26,VW13,ZZ19}
and the Narayana polynomials of types $A$ and $B$~\cite{YZ25}.
These results naturally raise the problem of whether the Delannoy polynomials admit a real stable multivariate refinement.
In this paper, we construct stable multivariate Delannoy polynomials in terms of weighted lattice paths.

Let $\cP_n$ be the set of Delannoy paths of size $n$.
For each $p\in\cP_n$,
assign a weight to each step of $p$,
and define $\wt(p)$ as the product of its step weights,
with the empty product equal to $1$.
If $\wt(H)=1$ and $\wt(V)=\wt(D)=x$,
then the weight sum of the paths in $\cP_n$ is the univariate Delannoy polynomial $d_n(x)$:
\[
d_n(x)=\sum_{k=0}^n d(n,k)x^k=\sum_{p\in\cP_n} \wt(p)
\]
(see \cite[Theorem 2.3]{WZC19} for example).
We now consider a multivariate refinement of $\wt(p)$
by distinguishing the steps of $p$ according to the antidiagonal lines on which they end.
For $0\le j\le n$, let
\[
L_j=\{(x,y)\mid x+y=j\}.
\]
Denote by $H_j$, $V_j$, and $D_j$ the horizontal, vertical, and diagonal steps, respectively,
that end on $L_j$.

\begin{defi}\label{defi}
	The $n$-th \defn{multivariate Delannoy polynomial} is defined as
	\[
	d_n(\bfx)=d_n(x_1,\ldots,x_n):=\sum_{p\in\cP_n}\wt(p),
	\]
	where the weights of steps in $p$ are given by
	\[
	\wt(H_j)=1,\quad \wt(V_j)=x_j, \quad
	\wt(D_j)=
	\begin{cases}
		x_{j-1}, & j\ \text{is even};\\
		x_j,     & j\ \text{is odd}.
	\end{cases}
	\]
\end{defi}

\begin{figure}[htbp]
	\centering
	\begin{tikzpicture}[x=0.45cm,y=0.45cm]
		\tikzset{
			ppath/.style={line width=0.9pt, line cap=round, line join=round},
			vtx/.style={circle, draw, fill=black, inner sep=1.2pt}
		}		
		\newcommand{\MiniPath}[5]{%
			\begin{scope}[shift={(#1,#2)}]
				\draw[ppath] #3;
				\foreach \P in {#4}{
					\node[vtx] at \P {};
				}
				\node at (1.25,-1.05) {$#5$};
			\end{scope}
		}		
		\MiniPath{0}{5}
		{(0,0)--(1,0)--(2,0)--(3,0)}
		{(0,0),(1,0),(2,0),(3,0)}
		{1}
		
		\MiniPath{4.5}{5}
		{(0,0)--(0,1)--(1,1)--(2,1)}
		{(0,0),(0,1),(1,1),(2,1)}
		{x_1}
		
		\MiniPath{9}{5}
		{(0,0)--(1,0)--(1,1)--(2,1)}
		{(0,0),(1,0),(1,1),(2,1)}
		{x_2}
		
		\MiniPath{13.5}{5}
		{(0,0)--(1,0)--(2,0)--(2,1)}
		{(0,0),(1,0),(2,0),(2,1)}
		{x_3}
		
		\MiniPath{18}{5}
		{(0,0)--(1,1)--(2,1)}
		{(0,0),(1,1),(2,1)}
		{x_1}
		
		\MiniPath{22.5}{5}
		{(0,0)--(1,0)--(2,1)}
		{(0,0),(1,0),(2,1)}
		{x_3}
		
		\MiniPath{0}{0}
		{(0,0)--(0,1)--(0,2)--(1,2)}
		{(0,0),(0,1),(0,2),(1,2)}
		{x_1x_2}
		
		\MiniPath{4.5}{0}
		{(0,0)--(0,1)--(1,1)--(1,2)}
		{(0,0),(0,1),(1,1),(1,2)}
		{x_1x_3}
		
		\MiniPath{9}{0}
		{(0,0)--(1,0)--(1,1)--(1,2)}
		{(0,0),(1,0),(1,1),(1,2)}
		{x_2x_3}
		
		\MiniPath{13.5}{0}
		{(0,0)--(1,1)--(1,2)}
		{(0,0),(1,1),(1,2)}
		{x_1x_3}
		
		\MiniPath{18}{0}
		{(0,0)--(0,1)--(1,2)}
		{(0,0),(0,1),(1,2)}
		{x_1x_3}
		
		\MiniPath{22.5}{0}
		{(0,0)--(0,1)--(0,2)--(0,3)}
		{(0,0),(0,1),(0,2),(0,3)}
		{x_1x_2x_3}		
	\end{tikzpicture}
	\caption{The weighted Delannoy paths in $\cP_3$.}
	\label{fig:P3}
\end{figure}
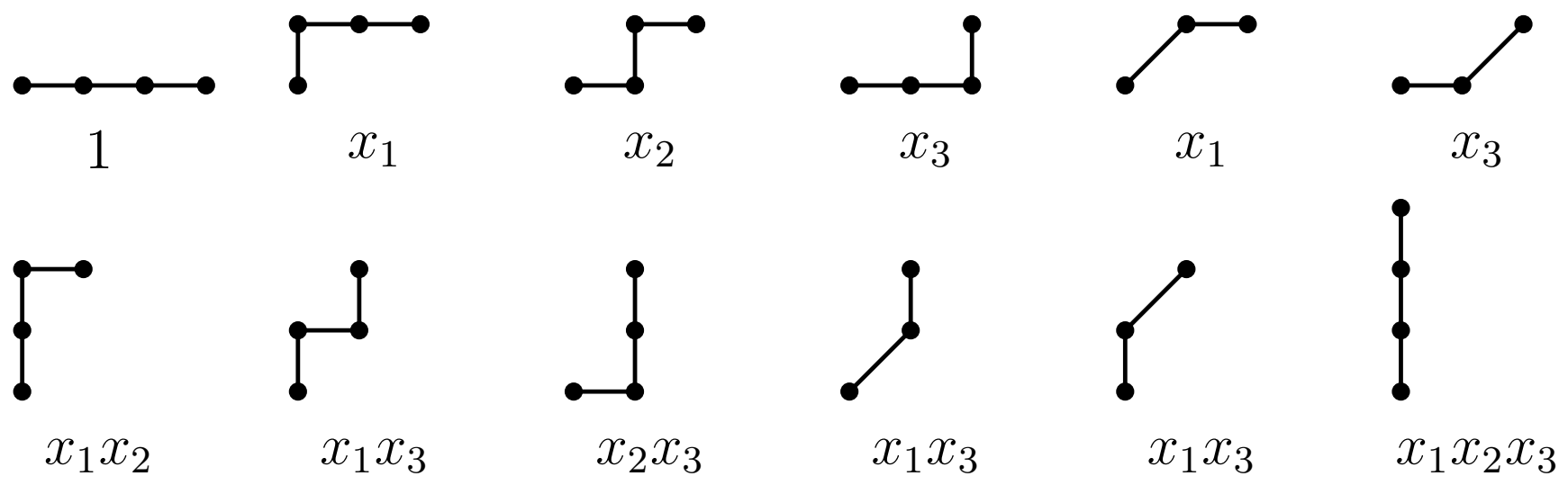

For example, the weighted Delannoy paths in $\cP_3$ are shown in Figure~\ref{fig:P3}.
Thus,
\[
d_3(\bfx)=1+2x_1+x_2+2x_3+x_1x_2+3x_1x_3+x_2x_3+x_1x_2x_3 .
\]

Under the diagonal specialization $x_1=\dots=x_n=x$,
these polynomials reduce to the univariate Delannoy polynomials.
The following theorem extends the real-rootedness of the univariate Delannoy polynomial to multivariate real stability.

\begin{thm}\label{thm:stability}
	For $n\geq 0$, the multivariate Delannoy polynomial $d_n(\bfx)$ is real stable.
\end{thm}

A polynomial is \defn{multi-affine} 
if its degree in each variable is at most one.
The weighting in Definition \ref{defi} makes $d_n(\bfx)$ multi-affine.
Indeed, if $j$ is even, $x_j$ can occur only as the weight of $V_j$;
if $j$ is odd, it can occur only as the weight of $V_j$, $D_j$, or $D_{j+1}$,
and no path can contain more than one of these steps.

Following~\cite{CW06},
the \defn{Rayleigh difference} of a multi-affine polynomial $f\in\R[x_1,\dots,x_n]$ with respect to $x_i$ and $x_j$ is defined as
\[
\Del_{ij}(f)=\frac{\partial f}{\partial x_i}\frac{\partial f}{\partial x_j}-f\frac{\partial^2 f}{\partial x_i\partial x_j}.
\]
Brändén~\cite{Bra07} characterized real stability of multi-affine polynomials through the nonnegativity of their Rayleigh differences.

\begin{thm}[{\cite[Theorem 5.6]{Bra07}}]
	\label{thm-Rayleigh}
	Let $f\in\R[x_1,\dots,x_n]$ be multi-affine.
	Then $f$ is real stable if and only if $\Del_{ij}(f)(x_1,\dots,x_n)\ge 0$ for all $(x_1,\dots,x_n)\in\R^n$ and $1\le i<j\le n$.
\end{thm}

Thus, Theorem~\ref{thm:stability} is equivalent to the nonnegativity of all Rayleigh differences of $d_n(\bfx)$ on $\R^n$. 
In fact, we prove the stronger statement that each of these Rayleigh differences is a perfect square.

\begin{thm}\label{thm:rayleigh-square}
	For $n\geq 2$ and $1\leq i<j\leq n$, we have 
	\[
	\Del_{ij}(d_n(\bfx))=d_{i-1}^2(\bfx) R_{jn}^2 \prod_{\substack{i<r<j\\ r\text{ odd}}}x_r^2
	\]
	where 
	\[
	R_{jn}=
	\begin{cases}
		d_{n-j}(x_{j+1},\ldots,x_n), &\text{$j$ is even};\\
		d_{n-j+1}(0,x_{j+1},\ldots,x_n), &\text{$j$ is odd}.
	\end{cases}
	\]
\end{thm}

The rest of the paper is organized as follows. 
Section \ref{sec:stability} gives an algebraic proof of Theorem~\ref{thm:stability} using the notion of proper position.
Section~\ref{sec:rayleigh} proves Theorem~\ref{thm:rayleigh-square} combinatorially by a weight-preserving injection on pairs of Delannoy paths.
Together with Theorem \ref{thm-Rayleigh}, this also proves Theorem~\ref{thm:stability}.

\section{Algebraic proof of Theorem \ref{thm:stability}}
\label{sec:stability}

To prove Theorem \ref{thm:stability}, we need the notion of proper position and a criterion for real stability via three-term recurrences.

Denote by $\cH_n(\C)$ (resp. $\cH_n(\R)$) the set of stable (resp. real stable) polynomials in $n$ variables.
Let $f,g\in\R[x_1,\dots,x_n]$.
We say that $f$ and $g$ are in \defn{proper position}, denoted by $f\ll g$, if $g+if\in\cH_n(\C)$.
A characterization of this notion was given by Borcea and Br\"and\'en~\cite{BB09}.

\begin{lem}[{\cite[Lemma 2.8]{BB09}}]
	\label{lem:BB-proper}
	Let $f,g\in\R[x_1,\ldots,x_n]$ be nonzero polynomials that are not constant multiples of each other,  and $x_{n+1}$ be a new indeterminate. 
	The following are equivalent:
	\begin{enumerate}[\rm(i)]
		\item $g\ll f;$
		\item $f+x_{n+1}g\in\cH_{n+1}(\R);$
		\item $f,g\in\cH_n(\R)$ and $\Ima\left(\frac{f(\bfx)}{g(\bfx)}\right)>0$ whenever $\Ima(x_i)>0$ for $1\le i\le n$.
	\end{enumerate}
\end{lem}

The following special case of \cite[Theorem 2.1]{HLY23} provides a useful criterion for proving real stability of polynomial sequences.

\begin{lem}[{\cite[Theorem 2.1]{HLY23}}]
	\label{lem:HLY}
	Suppose that $F,f,g\in\R[x_1,\ldots,x_n]$ satisfy
	\[
	F(\bfx)=a(\bfx)f(\bfx)+b(\bfx)g(\bfx),
	\]
	where $a(\bfx), b(\bfx) \in\cH_n(\R)$.
	If
	\begin{enumerate}[\rm(i)]
		\item $f(\bfx),g(\bfx)\in\cH_n(\R)$ and $g(\bfx)\ll f(\bfx),$
		\item $a(\bfx)\neq0$ and $\Ima\left(\frac{b(\bfx)}{a(\bfx)}\right)>0$ whenever $\Ima(x_i)>0$ for $1\le i\le n,$
	\end{enumerate}
	then $F(\bfx)\in\cH_n(\R)$.
\end{lem}

\begin{rem}
	For nonzero polynomials $f$ and $g$ that are not constant multiples of each other,
	condition (i) in Lemma \ref{lem:HLY} can be simplified to $g(\bfx)\ll f(\bfx)$ by Lemma \ref{lem:BB-proper}.
\end{rem}

Decomposing a Delannoy path according to its final step gives the following three-term recurrence.

\begin{prop}\label{prop:rec}
	For $n\ge2$, we have
	\begin{equation}\label{eq-rec}
		d_n(\bfx)=(1+x_n)d_{n-1}(\bfx)+
		\begin{cases}
			x_nd_{n-2}(\bfx),& n\text{ is odd};\\
			x_{n-1}d_{n-2}(\bfx),& n\text{ is even},
		\end{cases}
	\end{equation}
	with $d_0(\bfx)=1$ and $d_1(\bfx)=1+x_1$.
\end{prop}

\begin{proof}
Consider the final step of a path in $\cP_n$.
If it is horizontal or vertical,
deleting it leaves a path in $\cP_{n-1}$;
these two possibilities together contribute $(1+x_n)d_{n-1}(\bfx)$.
If it is diagonal, deleting it leaves a path in $\cP_{n-2}$,
and its weight is $x_n$ for odd $n$ and $x_{n-1}$ for even $n$.
Summing these contributions proves \eqref{eq-rec}.
\end{proof}

We are now ready to prove Theorem~\ref{thm:stability}.

\begin{proof}[Proof of Theorem~\ref{thm:stability}]
The constant polynomial $d_0(\bfx)=1$ is real stable.
To show that $d_n(\bfx)\in\cH_n(\R)$ for all $n\ge1$,
by Lemma \ref{lem:BB-proper} it suffices to prove the stronger statement that $d_{2m-1}(\bfx)\ll d_{2m}(\bfx)$ for all $m\ge1$,
as these polynomials are nonzero and are not constant multiples of each other.
We proceed by induction on $m$.

For $m=1$, we have $d_1(\bfx)=1+x_1$ and
\[
d_2(\bfx)=x_1+(1+x_1)(1+x_2).
\]
By definition, $d_1(\bfx)$ is real stable.
The Rayleigh difference $\Del_{12}(d_2(\bfx))=1$;
hence $d_2(\bfx)$ is real stable by Theorem \ref{thm-Rayleigh}.
Moreover,
\[
\Ima\left(\frac{d_2(\bfx)}{d_1(\bfx)}\right)
=\Ima(1+x_2)+\Ima\left(\frac{x_1}{1+x_1}\right)
=\Ima(x_2)+\frac{\Ima(x_1)}{|1+x_1|^2}>0
\]
whenever $\Ima(x_1), \Ima(x_2)>0.$
Thus, by Lemma~\ref{lem:BB-proper}, we have $d_1(\bfx)\ll d_2(\bfx) $.

Now suppose that $m\ge2$.
Assume that $d_{2m-3}(\bfx)\ll d_{2m-2}(\bfx)$.
We need to show that $d_{2m-1}(\bfx)\ll d_{2m}(\bfx)$,
or equivalently, by Lemma \ref{lem:BB-proper},
\begin{equation}\label{eq-F}
	F_m(\bfx)=F_m(x_1,\dots,x_{2m+1})
	:=d_{2m}(\bfx)+x_{2m+1} d_{2m-1}(\bfx)
\end{equation}
is real stable.
Rewriting \eqref{eq-rec} according to the parity of $n$,
we have
\begin{equation}\label{eq-rec-odd}
	d_{2m-1}(\bfx)=(1+x_{2m-1})d_{2m-2}(\bfx)+x_{2m-1}d_{2m-3}(\bfx)
\end{equation}
and
\begin{equation}\label{eq-rec-even}
	d_{2m}(\bfx)=(1+x_{2m})d_{2m-1}(\bfx)+x_{2m-1}d_{2m-2}(\bfx).
\end{equation}
Combining \eqref{eq-rec-odd} and \eqref{eq-rec-even} yields
\begin{equation}
	\label{eq-rec-even2}
	d_{2m}(\bfx)=[(1+x_{2m-1})(1+x_{2m})+x_{2m-1}]d_{2m-2}(\bfx)+x_{2m-1}(1+x_{2m})d_{2m-3}(\bfx).
\end{equation}
Substituting \eqref{eq-rec-odd} and \eqref{eq-rec-even2} into \eqref{eq-F} gives
\[
F_m(\bfx)=[u_m+(1+u_m)(1+v_m+w_m)]d_{2m-2}(\bfx)+u_m(1+v_m+w_m)d_{2m-3}(\bfx),
\]
where $u_m=x_{2m-1}$, $v_m=x_{2m}$ and $w_m=x_{2m+1}$.
Now apply Lemma \ref{lem:HLY} to $F_m(\bfx)$.
Since $d_2(\bfx)$ is real stable,
\[
u_m+(1+u_m)(1+v_m+w_m)=d_2(u_m,v_m+w_m)\in\cH_{2m+1}(\R).
\]
By definition, $u_m(1+v_m+w_m)\in\cH_{2m+1}(\R)$.
Moreover,
\begin{align*}
	\Ima\left(\frac{u_m+(1+u_m)(1+v_m+w_m)}{u_m(1+v_m+w_m)}\right)
	&=\Ima\left(\frac{1}{1+v_m+w_m}\right)+\Ima\left(1+\frac{1}{u_m}\right)\\
	&=-\frac{\Ima(v_m+w_m)}{|1+v_m+w_m|^2}-\frac{\Ima u_m}{|u_m|^2}
	<0,
\end{align*}
which implies
\[
\Ima\left(\frac{u_m(1+v_m+w_m)}{u_m+(1+u_m)(1+v_m+w_m)}\right)>0
\]
whenever $\Ima(u_m)$, $\Ima(v_m)$, $\Ima(w_m)>0$.
By the induction hypothesis, $d_{2m-3}(\bfx)\ll d_{2m-2}(\bfx)$.
Therefore, by Lemma \ref{lem:HLY}, $F_m(\bfx)\in\cH_{2m+1}(\R)$ as desired.
\end{proof}

\begin{rem}
In the proof of Theorem \ref{thm:stability},
we obtain $d_{2m-1}(\bfx)\ll d_{2m}(\bfx)$ for all $m\ge1$.
However, the relation $d_{2m}(\bfx)\ll d_{2m+1}(\bfx)$ does not hold in general.
Indeed, it fails when $m=1$.
In this case, we have $d_2(\bfx)=1+2x_1+x_2+x_1x_2$ and
\[
d_3(\bfx)=1+2x_1+x_2+2x_3+x_1x_2+3x_1x_3+x_2x_3+x_1x_2x_3.
\]
So
\begin{align*}
	d_3(\bfx)+x_4d_2(\bfx)
	=&1+2x_1+x_2+2x_3+x_4+x_1x_2+3x_1x_3+2x_1x_4+x_2x_3+x_2x_4\\
	&+x_1x_2x_3+x_1x_2x_4.
\end{align*}
A direct computation gives $\Del_{13}(d_3(\bfx)+x_4d_2(\bfx))=1+x_4$,
which is not always nonnegative.
So $d_3(\bfx)+x_4d_2(\bfx)$ is not real stable,
and hence $d_2\not\ll d_3$.
\end{rem}

We conclude this section with two remarks on the weights of diagonal steps.

\begin{rem}
\label{rem:naive-refinement}
A natural choice is to assign the same weight to vertical and diagonal steps ending on the same line $L_j$: 
\[
\wt(H_j)=1,\quad \wt(V_j)=\wt(D_j)=x_j.
\]
However, this weighting does not always produce real stable polynomials.  
For $n=3$, it gives 
\[
\what d_3(\bfx)=1+x_1+2x_2+2x_3+x_1x_2+2x_1x_3+2x_2x_3+x_1x_2x_3.
\]
Since $\Del_{13}(\what d_3(\bfx))=-x_2$ is not always nonnegative,
$\what d_3(\bfx)$ is not real stable.
Unlike a horizontal or vertical step, a diagonal step ending on $L_j$ must skip $L_{j-1}$.
This fact motivates us to treat the diagonal steps separately, as in Definition~\ref{defi}.
\end{rem}

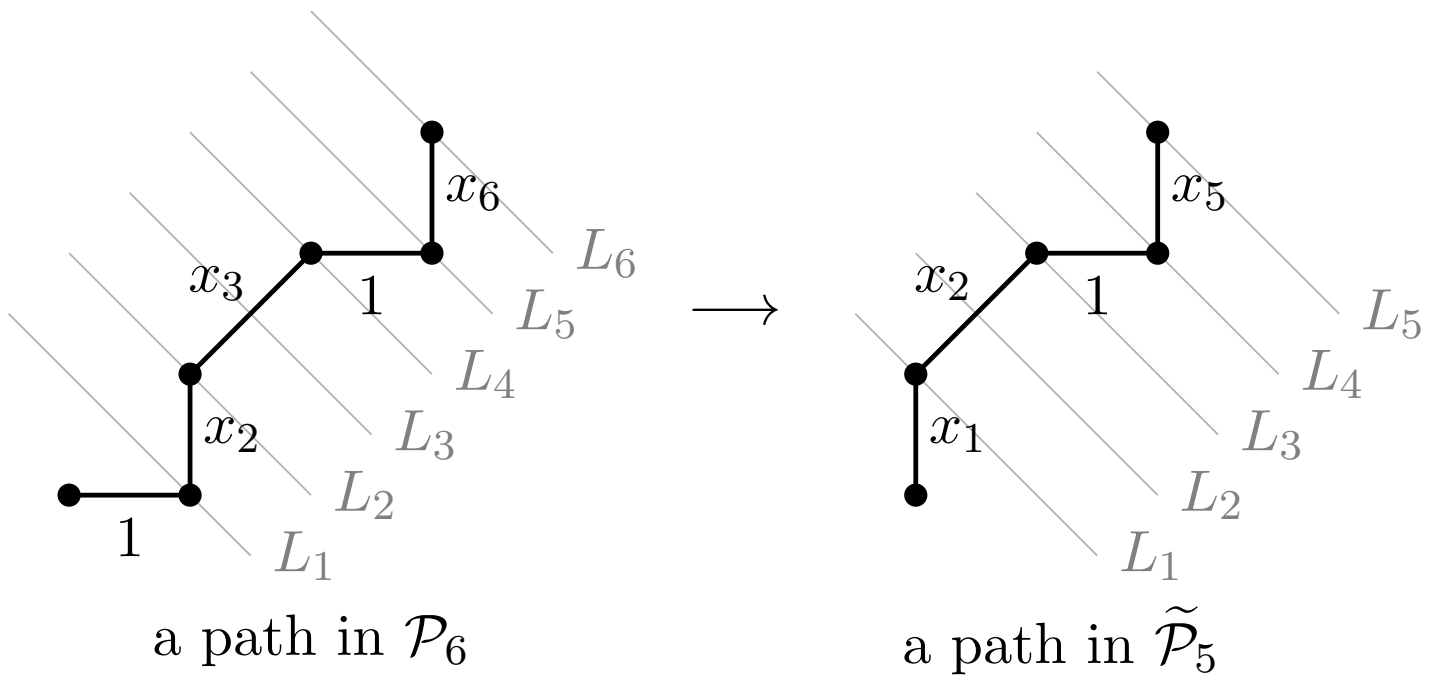
\begin{figure}[htbp]
	\centering 
	\vspace{-0.75em}
	\begin{tikzpicture}[
		scale=0.8,
		every node/.style={font=\small},
		adline/.style={gray!60,line width=0.35pt},
		pt/.style={circle, draw, fill=black, inner sep=1.4pt},
		lab/.style={font=\small},
		llab/.style={gray,font=\small},
		>=Latex
		]
		
		\begin{scope}[xshift=0cm, yshift=0cm]
			\foreach \x in {1,...,6}
				{
					\draw[adline] (0.5*\x-1,0.5*\x+1)--(0.5*\x+1,0.5*\x-1) node[llab,right] {$L_{\x}$};
				}
			\coordinate (A0) at (0,0);
			\coordinate (A1) at (1,0);
			\coordinate (A2) at (1,1);
			\coordinate (A3) at (2,2);
			\coordinate (A4) at (3,2);
			\coordinate (A5) at (3,3);
			
			\draw[line width=0.9pt] (A0)--(A1);
			\draw[line width=0.9pt] (A1)--(A2);
			\draw[line width=0.9pt] (A2)--(A3);
			\draw[line width=0.9pt] (A3)--(A4);
			\draw[line width=0.9pt] (A4)--(A5);
			
			\node[pt] at (A0) {};
			\node[pt] at (A1) {};
			\node[pt] at (A2) {};
			\node[pt] at (A3) {};
			\node[pt] at (A4) {};
			\node[pt] at (A5) {};
			
			\node[lab] at ($(A0)!0.5!(A1)+(0,-0.35)$) {$1$};
			\node[lab] at ($(A1)!0.5!(A2)+(0.35,0)$) {$x_2$};
			\node[lab,left] at ($(A2)!0.5!(A3)+(0.15,0.25)$) {$x_3$};
			\node[lab] at ($(A3)!0.5!(A4)+(0,-0.35)$) {$1$};
			\node[lab] at ($(A4)!0.5!(A5)+(0.35,0)$) {$x_6$};
			
			\node at (2,-1.2) {a path in $\cP_6$};
		\end{scope}
		
		\node at (5.5,1.5) {$\longrightarrow$};
		
		\begin{scope}[xshift=7cm, yshift=0cm]
			\foreach \x in {1,...,5}
				{
					\draw[adline] (0.5*\x-1,0.5*\x+1)--(0.5*\x+1,0.5*\x-1) node[llab,right] {$L_{\x}$};
				}
			\coordinate (B0) at (0,0);
			\coordinate (B1) at (0,1);
			\coordinate (B2) at (1,2);
			\coordinate (B3) at (2,2);
			\coordinate (B4) at (2,3);
			
			\draw[line width=0.9pt] (B0)--(B1);
			\draw[line width=0.9pt] (B1)--(B2);
			\draw[line width=0.9pt] (B2)--(B3);
			\draw[line width=0.9pt] (B3)--(B4);
			
			\node[pt] at (B0) {};
			\node[pt] at (B1) {};
			\node[pt] at (B2) {};
			\node[pt] at (B3) {};
			\node[pt] at (B4) {};
			
			\node[lab] at ($(B0)!0.5!(B1)+(0.35,0)$) {$x_1$};
			\node[lab,left] at ($(B1)!0.5!(B2)+(0.15,0.25)$) {$x_2$};
			\node[lab] at ($(B2)!0.5!(B3)+(0,-0.35)$) {$1$};
			\node[lab] at ($(B3)!0.5!(B4)+(0.35,0)$) {$x_5$};
			
			\node at (1.2,-1.2) {a path in $\wtil\cP_5$};		
		\end{scope}	
	\end{tikzpicture}
	\vspace{-0.7em}
	\caption{An illustration for $n=5$.}
	\label{fig:rem-bijection}
\end{figure}

\begin{rem}\label{rem:even-index}
In Definition~\ref{defi}, each diagonal step $D_j$ receives the odd-indexed variable among $x_{j-1}$ and $x_j$.
Choosing the even-indexed variable instead also gives real stable polynomials.
Let
\[
\wtil d_n(\bfx)=\sum_{p\in\wtil\cP_n}\wt(p),
\]
where $\wtil\cP_n$ is the set of Delannoy paths of size $n$ with step weights
\begin{equation}\label{eq-wt-even}
	\wt(H_j)=1,\quad \wt(V_j)=x_j, \quad
	\wt(D_j)=
	\begin{cases}
		x_j, & j\ \text{is even};\\
		x_{j-1},     & j\ \text{is odd}. 
	\end{cases}
\end{equation}
It follows that
\begin{equation}\label{eq-d&dtil}
	\wtil d_n(x_1,\ldots,x_n)=d_{n+1}(0,x_1,\ldots,x_n).
\end{equation}
Indeed, after the first variable is set to zero,
only paths beginning with $H_1$ have nonzero weight.
Deleting this initial step and translating the remaining path one unit to the left decreases every weight index by one,
giving \eqref{eq-d&dtil}.
Figure~\ref{fig:rem-bijection} illustrates this correspondence.
It is known~\cite{BB09} that specialization at a real constant preserves real stability.
Therefore, Theorem~\ref{thm:stability} implies that $\wtil d_n(\bfx)\in\cH_n(\R)$.
\end{rem}

\section{Combinatorial proof of Theorem \ref{thm:rayleigh-square}}
\label{sec:rayleigh}

We begin with an alternative expression of the Rayleigh differences.
Let $f\in\R[x_1,\dots,x_n]$ be multi-affine, and fix $1\leq i<j\leq n$. 
Then $f$ has a unique decomposition
\[
f=f_{ij}^{00}+x_i f_{ij}^{10}+x_j f_{ij}^{01}+x_i x_j f_{ij}^{11},
\]
where $f_{ij}^{st}\in\R[x_1,\ldots,\what{x_i},\ldots,\what{x_j},\ldots,x_n]$, for $s,t\in\{0,1\}$,
are independent of $x_i$ and $x_j$.
Direct calculation gives
\begin{equation}\label{eq-Del}
	\Del_{ij}(f)
	=\frac{\partial f}{\partial x_i}\frac{\partial f}{\partial x_j}-f\frac{\partial^2 f}{\partial x_i\partial x_j}
	=f_{ij}^{01}f_{ij}^{10}-f_{ij}^{00}f_{ij}^{11}.
\end{equation}

Applying this decomposition to $d_n(\bfx)$, we write
\begin{equation}\label{eq-dn-decom}
	d_n(\bfx)=d_{ij}^{00}(\bfx)+x_i d_{ij}^{10}(\bfx)+x_j d_{ij}^{01}(\bfx)+x_i x_j d_{ij}^{11}(\bfx).
\end{equation}
For example, when $n=3$, $i=1$, and $j=2$, we have
\[
d_3(\bfx)=(1+2x_3)+x_1(2+3x_3)+x_2(1+x_3)+x_1x_2(1+x_3).
\]
Here $d_{12}^{00}=1+2x_3$, $d_{12}^{10}=2+3x_3$, $d_{12}^{01}=1+x_3$, $d_{12}^{11}=1+x_3$.
Thus,
\[
\Del_{12}(d_3(\bfx))=(1+x_3)(2+3x_3)-(1+x_3)(1+2x_3)=(1+x_3)^2.
\]

For $s,t\in\{0,1\}$, let $\cP_{ij}^{st}$ be the set of Delannoy paths $p\in\cP_n$ such that $x_i\mid\wt(p)$
if and only if $s=1$, and
$x_j\mid\wt(p)$
if and only if $t=1$.
Then \eqref{eq-dn-decom} implies that
\[
\begin{aligned}
	d_{ij}^{00}(\bfx)&=\sum_{p\in\cP_{ij}^{00}}\wt(p), \quad & x_i d_{ij}^{10}(\bfx)&=\sum_{p\in\cP_{ij}^{10}}\wt(p),\\
	x_j d_{ij}^{01}(\bfx)&=\sum_{p\in\cP_{ij}^{01}}\wt(p), \quad & x_i x_j d_{ij}^{11}(\bfx)&=\sum_{p\in\cP_{ij}^{11}}\wt(p).
\end{aligned}
\]
Combining this with \eqref{eq-Del}, we obtain
\[
x_i x_j\Del_{ij}(d_n(\bfx))
=\sum_{p\in\cP_{ij}^{01}}\wt(p)
\sum_{q\in\cP_{ij}^{10}}\wt(q)
-\sum_{p\in\cP_{ij}^{00}}\wt(p)
\sum_{q\in\cP_{ij}^{11}}\wt(q),
\]
or equivalently,
\begin{equation}
	\label{eq:rayleigh-path-pairs}
	x_ix_j\Del_{ij}(d_n(\bfx))=
	\sum_{(p,q)\in\cP_{ij}^{01}\times\cP_{ij}^{10}}\wt(p)\wt(q)-\sum_{(p,q)\in\cP_{ij}^{00}\times\cP_{ij}^{11}}\wt(p)\wt(q).
\end{equation}
Thus $x_i x_j\Del_{ij}(d_n(\bfx))$ is the signed weight sum of pairs of Delannoy paths.
We now construct a weight-preserving injection
\[
\psi:\cP_{ij}^{00}\times\cP_{ij}^{11}\to\cP_{ij}^{01}\times\cP_{ij}^{10}
\]
and show that the total weight of the pairs in $(\cP_{ij}^{01}\times\cP_{ij}^{10})\setminus\psi(\cP_{ij}^{00}\times\cP_{ij}^{11})$ is $x_ix_j$ times a perfect square.

Recall that $L_k=\{(x,y)\mid x+y=k\}$.
We say that a path $p$ \defn{meets} $L_k$ if $p$ has a vertex on $L_k$, and $p$ \defn{skips} $L_k$ otherwise. 
Note that a path skips $L_k$ precisely when it reaches $L_{k+1}$ by $D_{k+1}$.
The following proposition ensures that the map $\psi$ can be defined.

\begin{prop}\label{prop-commonLk}
	For each $(p,q)\in\cP_{ij}^{00}\times\cP_{ij}^{11}$, 
	there exists an integer $k$ with $i\leq k<j$ such that both $p$ and $q$ meet $L_k$.
\end{prop}

\begin{proof}
There are four parity cases of $(i,j)$.
We give the details only for the case in which both $i$ and $j$ are odd;
the other cases are treated similarly and summarized in Table~\ref{tab:parity-cases}.

Suppose, for a contradiction, that there is no integer $k$ with $i\le k<j$ such that both $p$ and $q$ meet $L_k$.
Note that $p\in\cP_{ij}^{00}$ has neither $x_i$ nor $x_j$ in the step weights.
Since $i$ is odd, $p$ must meet $L_i$ via $H_i$:
the alternatives $V_i$, $D_i$, and skipping $L_i$ via $D_{i+1}$ all contribute $x_i$.
Hence $q$ skips $L_i$ and meets $L_{i+1}$ via $D_{i+1}$,
forcing $p$ to skip $L_{i+1}$.
A path cannot skip two consecutive lines.
Thus, whenever one path skips a line, it must meet the next one, forcing the other path to skip that next line.
Consequently, $p$ and $q$ meet $L_i,L_{i+1},\ldots,L_{j-1}$ alternately.
As both $i$ and $j$ are odd, $p$ skips $L_{j-1}$ and reaches $L_j$ via $D_j$.
This step has weight $x_j$, contradicting the choice of $p$.
\end{proof}

\begin{table}[htbp]
	\centering
  \caption{The parity-dependent contradictions}
	\begin{tabular}{@{}l l|c|c@{}}
		\toprule
		\multicolumn{2}{c|}{$(i,j)$ parity}
		& Initial behavior at $L_i$
		& Contradiction \\
		\midrule
		(1) & $i$ odd, $j$ odd
		& $p$ meets $L_i$ via $H_i$
		& $p$ meets $L_j$ via $D_j$ $\implies x_j\mid\wt(p)$\\
		(2) & $i$ odd, $j$ even
		& $p$ meets $L_i$ via $H_i$
		& $q$ meets $L_j$ via $D_j$ $\implies x_j\nmid\wt(q)$ \\
		(3) & $i$ even, $j$ odd
		& $q$ meets $L_i$ via $V_i$
		& $p$ meets $L_j$ via $D_j$ $\implies x_j\mid\wt(p)$ \\
		(4) & $i$ even, $j$ even
		& $q$ meets $L_i$ via $V_i$
		& $q$ meets $L_j$ via $D_j$ $\implies x_j\nmid\wt(q)$ \\
		\bottomrule
	\end{tabular}
	\label{tab:parity-cases}
\end{table}

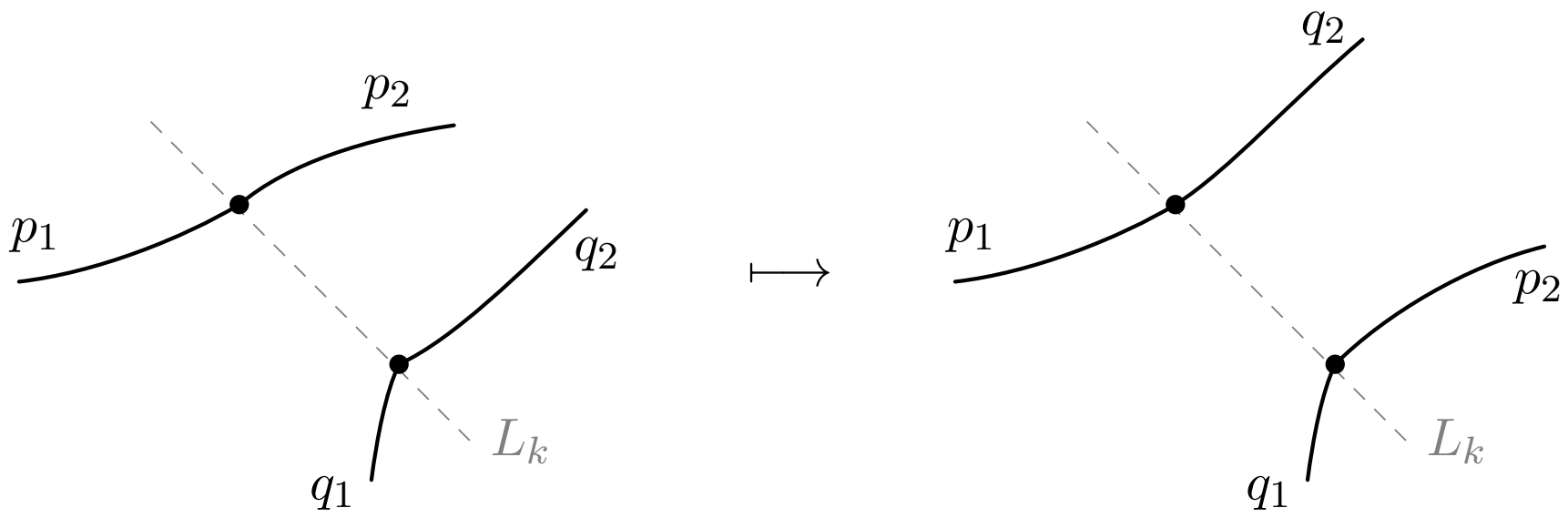
\begin{figure}[htbp]
	\centering	
	\begin{tikzpicture}[scale=0.8,
		x=1cm,
		y=1cm,
		>=Latex,
		line cap=round,
		line join=round,
		every node/.style={font=\small},
		adline/.style={dashed,gray,line width=0.35pt},
    pathnum/.style={black},
		]
		
		\begin{scope}[xshift=0cm]
			
			\draw[adline] (1.20,2.75) -- (4.10,-0.15) node[right] {$L_k$};
			
			\coordinate (A) at (2.00,2.00);
			\coordinate (B) at (3.45,0.55);
			
			\draw[black,thick]
			(0.00,1.30)
			.. controls (0.70,1.38) and (1.50,1.70) ..
			(A)
			.. controls (2.45,2.38) and (3.15,2.60) ..
			(3.95,2.72);
			
			\draw[black,thick]
			(3.20,-0.50)
			.. controls (3.25,-0.12) and (3.33,0.30) ..
			(B)
			.. controls (3.95,0.78) and (4.55,1.38) ..
			(5.15,1.95);
			
			\fill[black] (A) circle (2.5pt);
			\fill[black] (B) circle (2.5pt);
			
			\node[pathnum] at (0.15,1.72) {$p_1$};
			\node[pathnum] at (3.35,3.02) {$p_2$};
			\node[pathnum] at (2.85,-0.62) {$q_1$};
			\node[pathnum] at (5.25,1.55) {$q_2$};
			
			\node[pathnum] at (2.55,-1.30)
			{$(p,q)\in\cP_{ij}^{00}\times\cP_{ij}^{11}$};	
		\end{scope}
			
		\node[pathnum] at (7,1.35) {$\longmapsto$};
				
		\begin{scope}[xshift=8.5cm]
			
			\draw[adline] (1.20,2.75) -- (4.10,-0.15) node[right] {$L_k$};
			
			\coordinate (A2) at (2.00,2.00);
			\coordinate (B2) at (3.45,0.55);
			
			\draw[black,thick]
			(0.00,1.30)
			.. controls (0.70,1.38) and (1.50,1.70) ..
			(A2)
			.. controls (2.50,2.32) and (3.05,2.95) ..
			(3.70,3.50);
			
			\draw[black,thick]
			(3.20,-0.50)
			.. controls (3.25,-0.12) and (3.33,0.30) ..
			(B2)
			.. controls (3.90,1.00) and (4.65,1.45) ..
			(5.35,1.62);
			
			\fill[black] (A2) circle (2.5pt);
			\fill[black] (B2) circle (2.5pt);
			
			\node[pathnum] at (0.15,1.72) {$p_1$};
			\node[pathnum] at (3.35,3.60) {$q_2$};
			\node[pathnum] at (2.85,-0.62) {$q_1$};
			\node[pathnum] at (5.30,1.25) {$p_2$};
			
			\node[pathnum] at (2.55,-1.30)
			{$(p',q')\in\cP_{ij}^{01}\times\cP_{ij}^{10}$};
		\end{scope}
	\end{tikzpicture}
	\caption{The map $\psi$ with $p=p_1p_2$, $q=q_1q_2$ and $p'=p_1q_2$, $q'=q_1p_2$.}
	\label{fig:switching}
\end{figure}

Choose the least $k$ in Proposition \ref{prop-commonLk}.
Cut the paths $p$ and $q$ at their vertices on $L_k$,
writing $p=p_1p_2$ and $q=q_1q_2$.
Switching the tail step sequences $p_2$ and $q_2$ and translating each tail to its new initial vertex defines a map
\[
\begin{array}{cccc}
	\psi: & \cP_{ij}^{00}\times\cP_{ij}^{11} & \longrightarrow & \cP_{ij}^{01}\times\cP_{ij}^{10}\\[5pt]
	& (p,q) & \longmapsto & (p',q'),
\end{array}
\]
where $p'=p_1q_2$ and $q'=q_1p_2$.
See Figure~\ref{fig:switching} for an illustration.

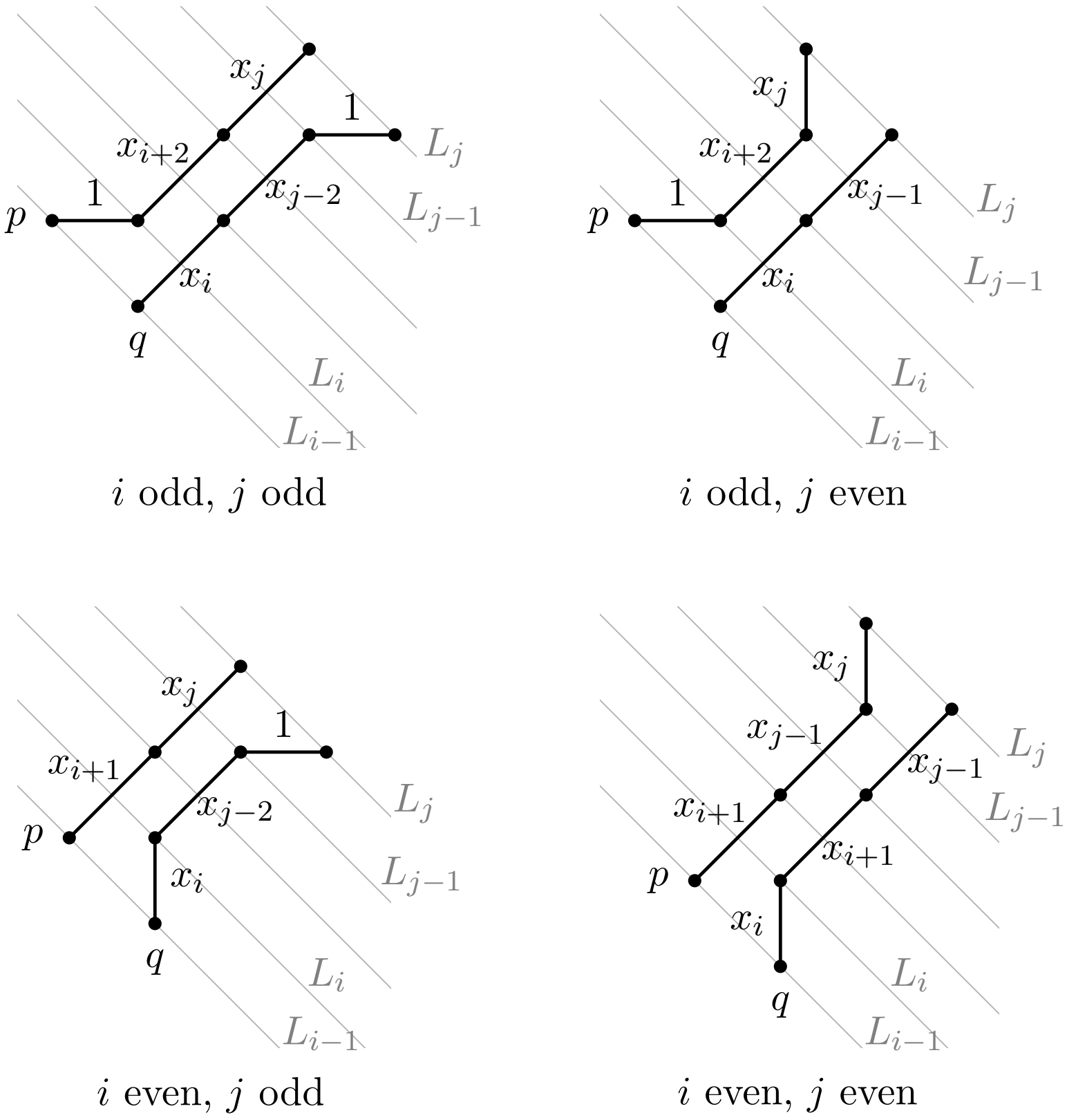
\begin{figure}[htbp]
	\centering	
	\begin{tikzpicture}[
		x=0.82cm,
		y=0.82cm,
		line cap=round,
		line join=round,
		pathline/.style={black,line width=0.9pt},
		adline/.style={gray!60,line width=0.35pt},
		wlab/.style={black,font=\small},
		llab/.style={gray,font=\small},
		every node/.style={font=\small}
		]
		
		\newcommand{\gridsix}{
			\begin{scope}
				\clip (-0.30,-2.10) rectangle (4.35,3.05);
				
				\foreach \c in {0.65,1.65,2.65,3.65,4.65,5.65}{
					\draw[adline]
					(-2,{\c+2})
					--
					(6,{\c-6});
				}
			\end{scope}
			
			\node[llab,right] at (4.25, 1.40) {$L_j$};
			\node[llab,right] at (4.00, 0.65) {$L_{j-1}$};
			
			\node[llab,right] at (2.90,-1.25) {$L_i$};
			\node[llab,right] at (2.60,-1.95) {$L_{i-1}$};
		}
		
		\newcommand{\gridfive}{
			\begin{scope}
				\clip (-0.30,-2.10) rectangle (4.05,3.05);
				
				\foreach \c in {0.65,1.65,2.65,3.65,4.65}{
					\draw[adline]
					(-2,{\c+2})
					--
					(6,{\c-6});
				}
			\end{scope}
			
			\node[llab,right] at (3.90, 0.75) {$L_j$};
			\node[llab,right] at (3.75,-0.10) {$L_{j-1}$};
			
			\node[llab,right] at (2.90,-1.25) {$L_i$};
			\node[llab,right] at (2.60,-1.95) {$L_{i-1}$};
		}
		
		\begin{scope}[shift={(0,0)}]
			
			\gridsix
			
			\coordinate (ooA) at (0.10,0.55);   
			\coordinate (ooB) at (1.10,0.55);   
			\coordinate (ooC) at (2.10,1.55);
			\coordinate (ooD) at (3.10,2.55);   
			
			\draw[pathline]
			(ooA)--(ooB)--(ooC)--(ooD);
			
			\coordinate (ooE) at (1.10,-0.45);  
			\coordinate (ooF) at (2.10,0.55);
			\coordinate (ooG) at (3.10,1.55);   
			\coordinate (ooH) at (4.10,1.55);   
			
			\draw[pathline]
			(ooE)--(ooF)--(ooG)--(ooH);
			
			\foreach \P in {ooA,ooB,ooC,ooD,ooE,ooF,ooG,ooH}{
				\fill (\P) circle (1.8pt);
			}
			
			\node[left=3pt]  at (ooA) {$p$};
			\node[below=3pt] at (ooE) {$q$};
			
			\node[wlab,above] at ($(ooA)!0.50!(ooB)$) {$1$};
			
			\node[wlab,left] at ($(ooB)!0.80!(ooC)$) {$x_{i+2}$};
			\node[wlab,left] at ($(ooC)!0.70!(ooD)$) {$x_j$};
			
			\node[wlab,right] at ($(ooE)!0.30!(ooF)$) {$x_i$};
			\node[wlab,right] at ($(ooF)!0.30!(ooG)$) {$x_{j-2}$};
			
			\node[wlab,above] at ($(ooG)!0.50!(ooH)$) {$1$};
			
			\node at (2.05,-2.65)
			{$i$ odd, $j$ odd};
			
		\end{scope}
		
		\begin{scope}[shift={(6.8,0)}]
			
			\gridfive
			
			\coordinate (oeA) at (0.10,0.55);   
			\coordinate (oeB) at (1.10,0.55);   
			\coordinate (oeC) at (2.10,1.55);   
			\coordinate (oeD) at (2.10,2.55);   
			
			\draw[pathline]
			(oeA)--(oeB)--(oeC)--(oeD);
			
			\coordinate (oeE) at (1.10,-0.45);  
			\coordinate (oeF) at (2.10,0.55);
			\coordinate (oeG) at (3.10,1.55);   
			
			\draw[pathline]
			(oeE)--(oeF)--(oeG);
			
			\foreach \P in {oeA,oeB,oeC,oeD,oeE,oeF,oeG}{
				\fill (\P) circle (1.8pt);
			}
			
			\node[left=3pt]  at (oeA) {$p$};
			\node[below=3pt] at (oeE) {$q$};
			
			\node[wlab,above] at ($(oeA)!0.50!(oeB)$) {$1$};
			
			\node[wlab,left] at ($(oeB)!0.80!(oeC)$) {$x_{i+2}$};
			\node[wlab,left] at ($(oeC)!0.50!(oeD)$) {$x_j$};
			
			\node[wlab,right] at ($(oeE)!0.30!(oeF)$) {$x_i$};
			\node[wlab,right] at ($(oeF)!0.30!(oeG)$) {$x_{j-1}$};
			
			\node at (1.95,-2.65)
			{$i$ odd, $j$ even};
			
		\end{scope}
		
		\begin{scope}[shift={(0,-7.0)}]
			
			\gridfive
			
			\begin{scope}[shift={(0.20,-0.20)}]

			\coordinate (eoA) at (0.10,0.55);   
			\coordinate (eoB) at (1.10,1.55);
			\coordinate (eoC) at (2.10,2.55);   
			
			\draw[pathline]
			(eoA)--(eoB)--(eoC);
			
			\coordinate (eoE) at (1.10,-0.45);  
			\coordinate (eoF) at (1.10,0.55);   
			\coordinate (eoG) at (2.10,1.55);   
			\coordinate (eoH) at (3.10,1.55);   
			
			\draw[pathline]
			(eoE)--(eoF)--(eoG)--(eoH);
			
			\foreach \P in {eoA,eoB,eoC,eoE,eoF,eoG,eoH}{
				\fill (\P) circle (1.8pt);
			}
			
			\node[left=3pt]  at (eoA) {$p$};
			\node[below=3pt] at (eoE) {$q$};
			
			\node[wlab,left] at ($(eoA)!0.80!(eoB)$) {$x_{i+1}$};
			\node[wlab,left] at ($(eoB)!0.70!(eoC)$) {$x_j$};
			
			\node[wlab,right] at ($(eoE)!0.50!(eoF)$) {$x_i$};
			\node[wlab,right] at ($(eoF)!0.30!(eoG)$) {$x_{j-2}$};
			\node[wlab,above] at ($(eoG)!0.50!(eoH)$) {$1$};
			\end{scope}
			\node at (1.95,-2.65)
			{$i$ even, $j$ odd};
			
		\end{scope}
		
		\begin{scope}[shift={(6.8,-7.0)}]
			
			\gridsix
			
			\begin{scope}[shift={(0.4,-0.4)}]
				
				\coordinate (eeA) at (0.40,0.25);
				\coordinate (eeB) at (1.40,1.25);
				\coordinate (eeC) at (2.40,2.25);
				\coordinate (eeD) at (2.40,3.25);
				
				\draw[pathline]
				(eeA)--(eeB)--(eeC)--(eeD);
				
				\coordinate (eeE) at (1.40,-0.75);
				\coordinate (eeF) at (1.40,0.25);
				\coordinate (eeG) at (2.40,1.25);
				\coordinate (eeH) at (3.40,2.25);
				
				\draw[pathline]
				(eeE)--(eeF)--(eeG)--(eeH);
				
				\foreach \P in {eeA,eeB,eeC,eeD,eeE,eeF,eeG,eeH}{
					\fill (\P) circle (1.8pt);
				}
				
				\node[left=3pt] at (eeA) {$p$};
				\node[below=3pt] at (eeE) {$q$};
				
				\node[wlab,left] at ($(eeA)!0.80!(eeB)$) {$x_{i+1}$};
				\node[wlab,left] at ($(eeB)!0.70!(eeC)$) {$x_{j-1}$};
				\node[wlab,left] at ($(eeC)!0.50!(eeD)$) {$x_j$};
				
				\node[wlab,left] at ($(eeE)!0.50!(eeF)$) {$x_i$};
				\node[wlab,right] at ($(eeF)!0.30!(eeG)$) {$x_{i+1}$};
				\node[wlab,right] at ($(eeG)!0.30!(eeH)$) {$x_{j-1}$};
				
			\end{scope}
			
			\node at (2.00,-2.65) {$i$ even, $j$ even};
			
		\end{scope}
		
	\end{tikzpicture}	
	\caption{Steps of $(p,q)$ between $L_{i-1}$ and $L_j$.}
	\label{fig:four-middle-configurations}
\end{figure}

Since each path has a vertex on the common cut $L_k$,
neither path contains a diagonal step that skips $L_k$.
So all occurrences of $x_i$ lie in the prefixes,
whereas all occurrences of $x_j$ lie in the tails.
Thus the exchange produces a pair in $\cP_{ij}^{01}\times\cP_{ij}^{10}$.
Translating a tail between two vertices on the same antidiagonal line does not change the antidiagonal line on which any of its steps ends.
Hence the exchange preserves the product of the weights:
\[
\wt(p)\wt(q)=\wt(p')\wt(q').
\]
The prefixes and the least common cut $L_k$ are unchanged,
so exchanging the tails again recovers the original pair.
Therefore, a pair is in
\[
(\cP_{ij}^{01}\times\cP_{ij}^{10})\setminus\psi(\cP_{ij}^{00}\times\cP_{ij}^{11})
\]
if and only if the two paths have no common meeting line $L_k$ with $i\le k<j$.
The two paths then meet the intermediate lines $L_i,L_{i+1},\ldots,L_{j-1}$ alternately.
Their middle portions are forced as in Figure~\ref{fig:four-middle-configurations},
including the horizontal or vertical boundary steps.

Let $\cQ_{ij}^{01}$ and $\cQ_{ij}^{10}$ consist of paths whose middle portions have the respective forms shown in Figure~\ref{fig:four-middle-configurations}.
Then
\[
\cQ_{ij}^{01}\times\cQ_{ij}^{10}=(\cP_{ij}^{01}\times\cP_{ij}^{10}) \setminus \psi(\cP_{ij}^{00}\times\cP_{ij}^{11}).
\]
Combining this identity with \eqref{eq:rayleigh-path-pairs}, we obtain
\begin{equation}\label{eq-xxDelta}
	x_i x_j\Del_{ij}(d_n(\bfx))
	=\sum_{(p,q)\in\cQ_{ij}^{01}\times\cQ_{ij}^{10}}\wt(p)\wt(q)
	=\sum_{p\in\cQ_{ij}^{01}}\wt(p)
	\sum_{q\in\cQ_{ij}^{10}}\wt(q).
\end{equation}

Each pair $(p,q)\in\cQ_{ij}^{01}\times\cQ_{ij}^{10}$ can be decomposed uniquely as $p=p_1p_2p_3$ and $q=q_1q_2q_3$,
where $p_1, q_1\in\cP_{i-1}$,
the middle parts $p_2, q_2$ are forced by Figure~\ref{fig:four-middle-configurations},
and $p_3, q_3$ are path segments from $L_j$ to $L_n$.
Clearly,
\[
\sum_{p_1} \wt(p_1)=\sum_{q_1} \wt(q_1)=d_{i-1}(\bfx).
\]
The forced middle portions have weights
\[
\wt(p_2)=x_j \prod_{\substack{i<r<j\\ r\text{ odd}}}x_r,
\qquad
\wt(q_2)=x_i \prod_{\substack{i<r<j\\ r\text{ odd}}}x_r.
\]
For each suffix $p_3$ or $q_3$, we translate the initial vertex on $L_j$ to the origin, and relabel $L_{j+r}$ as $L_r$ for $0\leq r\leq n-j.$
The resulting paths have size $n-j$,
with step weights:
\begin{equation}\label{eq-new-wt}
	\wt(H_r)=1,\quad
	\wt(V_r)=x_{j+r},\quad
	\wt(D_r)=
	\begin{cases}
		x_{j+r-1}, & j+r\text{ is even};\\
		x_{j+r}, & j+r\text{ is odd}.
	\end{cases}
\end{equation}
Note that in Definition \ref{defi}, the weight of each $D_r$ always uses the odd-indexed variable among $x_r$ and $x_{r-1}$.
For even $j$, the weighting in \eqref{eq-new-wt} gives precisely the original weights shifted by $j$.
Thus we have
\[
\sum_{p_3} \wt(p_3)=\sum_{q_3} \wt(q_3)=d_{n-j}(x_{j+1},\ldots,x_n).
\]
For odd $j$, the diagonal step $D_r$ has weight $x_{j+r}$ for even $r$ and $x_{j+r-1}$ for odd $r$.
This is the even-indexed version of the original weights (see Remark \ref{rem:even-index}) shifted by $j$.
Hence,
\[
\sum_{p_3} \wt(p_3)=\sum_{q_3} \wt(q_3)=\wtil d_{n-j}(x_{j+1},\ldots,x_n)=d_{n-j+1}(0,x_{j+1},\ldots,x_n).
\]

Note that the choices of $p_1$ and $p_3$ are independent,
and likewise for $q_1$ and $q_3$.
Hence,
\begin{align}\label{eq-wt-p}
	\nonumber
	\left(\sum_{p\in\cQ_{ij}^{01}}\wt(p)\right)
	=&\sum_{p=p_1p_2p_3} \wt(p_1)\wt(p_2)\wt(p_3)\\
	\nonumber
	=&x_j \prod_{\substack{i<r<j\\ r\text{ odd}}}x_r \sum_{p_1} \wt(p_1) \sum_{p_3} \wt(p_3)\\
	=&x_j d_{i-1}(\bfx) R_{jn} \prod_{\substack{i<r<j\\ r\text{ odd}}}x_r,
\end{align}
where
\[
R_{jn}=
\begin{cases}
	d_{n-j}(x_{j+1},\ldots,x_n), &\text{$j$ is even},\\
	d_{n-j+1}(0,x_{j+1},\ldots,x_n), &\text{$j$ is odd}.
\end{cases}
\]
Similarly,
\begin{align}\label{eq-wt-q}
	\sum_{q\in\cQ_{ij}^{10}}\wt(q)
	=x_i d_{i-1}(\bfx) R_{jn} \prod_{\substack{i<r<j\\ r\text{ odd}}}x_r.
\end{align}
Substituting \eqref{eq-wt-p} and \eqref{eq-wt-q} into \eqref{eq-xxDelta} and canceling the factor $x_i x_j$, we derive
\[
\Del_{ij}(d_n(\bfx))=d_{i-1}^2(\bfx) R_{jn}^2 \prod_{\substack{i<r<j\\ r\text{ odd}}}x_r^2.
\]
This completes the proof of Theorem~\ref{thm:rayleigh-square}.

\end{document}